\documentclass[a4paper]{article}
\usepackage{amsmath}
\usepackage{amsfonts}
\usepackage{amssymb}
\usepackage{amstext}
\usepackage{latexsym}

\newcommand{\whitebox}{$\Box$}

\newenvironment{proof}[1][Proof]{\textbf{#1.} }{\hfill\whitebox\bigskip}

\newcommand{\RR}{\mathbb{R}}

\newcommand{\NN}{\mathbb{N}}
\newcommand{\RRR}{\mathbb{R}^n}

\begin{document}
\title{Critical sets of nonlinear Sturm-Liouville operators of 
Ambrosetti-Prodi type}
\author{H. Bueno\quad and\quad C. Tomei\thanks{Supported in 
part by CNPq-Brazil, FINEP and FAPERJ.}}
\date{\small \textit{Dep. de Matem\'{a}tica - ICEx, UFMG, Caixa Postal 702, \linebreak Belo Horizonte, MG 30.123-970, Brazil.\linebreak \textit{Dep. de Matem\'atica, PUC-Rio, R. Marqu\^es de S. Vicente 225}\linebreak Rio de Janeiro, RJ 22453-900, Brazil.
\linebreak email - H. Bueno: hamilton@mat.ufmg.br\quad C. Tomei: tomei@mat.puc-rio.br}}
\maketitle

\begin{abstract}The critical set $C$ of the operator $F:H^2_D([0,\pi])\to L^2([0,\pi])$ defined by $F(u)=-u''+f(u)$ is studied. Here $X:=H^2_D([0,\pi])$ stands for the set of functions that satisfy the Dirichlet boundary conditions and whose derivatives are in $L^2([0,\pi])$. For generic
nonlinearities $f$, $C=\cup C_k$ decomposes into manifolds of codimension 1 in $X$. 
If $f''<0$ or $f''>0$, the set $C_j$ is shown to be non-empty if, and only if, $-j^2$ 
(the $j$-th eigenvalue of $u''$) is in the range of $f'$. The critical components $C_k$ are 
(topological) hyperplanes. 
\end{abstract}
\markright{H. Bueno and C. Tomei}

Let $\Omega\subset \RRR$ be a bounded, smooth domain, $f\in C^k(\RR)$, $k\geq 2$, and $g:\Omega\to \mathbb{R}$ be two given functions. Since the pioneering work of Ambrosetti-Prodi (\cite{AP}), the nonlinear problem
\[-\Delta u +f(u)=g,\quad u|_{\partial\Omega}=0,\]
has been studied by various methods of nonlinear analysis in the case
\[\sigma(-\Delta)\cap f'(\mathbb{R})\neq \emptyset.\]
($\sigma(-\Delta)$ stands for the spectrum of $-\Delta$). 

In this paper we deal with the one dimensional version of this problem, namely
\begin{equation}\label{dim1}
-u''+f(u)=g(t),\quad u(0)=u(\pi)=0.
\end{equation}
There exists an extensive literature concerning this problem (see, for example, \cite{C}, 
\cite{R}, and references therein). 
Usually one obtains {\it a priori} estimates on the nonlinearity, which  result in bounds on the number of solutions of $(\ref{dim1})$. 

Despite the success of the different methods used in this problem, they do not give insight into the nature of the change in the number of solutions of $(\ref{dim1})$, which was one of the major features of the work of Ambrosetti and Prodi, as well as subsequent work of 
Berger and Podolak (\cite{BP}). 
There is a certain reappraisal in recent years (\cite{BZ}, \cite{M},\cite{R}) of the Ambrosetti-Prodi method. This paper, whose approach is also geometric, 
may be seen as a portion of a larger project, inspired by some techniques and methods also 
present in \cite{MST}. For a first order differential equation, these authors characterized
the critical set $C$, then studied of the stratification of $C$ by
different Morin singularities and finally considered the geometry of the 
image of the critical set.  

Here we address only the first part of this project: for convex or concave nonlinearities, 
we characterize the critical set $C$ of the Sturm-Liouville operator
\[\begin{array}{cccc}F:&X&\to &Y\\ &u&\mapsto&-u''+f(u),
\end{array}\]
defined in the Sobolev space $X:=H^2_D([0,\pi])$ of functions that satisfy the Dirichlet boundary conditions and whose second derivatives are in $Y:=L^2([0,\pi])$. 
The results also hold for different pairs of spaces without difficulty.

Our main result may be synthesized as follows
\bigskip

\noindent {\bf Theorem A}: {\it Suppose that the nonlinearity $f$ satisfies $f''>0$ or $f''<0$. 
Then the critical set $C$ of the operator $F$ decomposes into connected components $C_j$, 
associated to the free eigenvalues $\{1^2, 2^2, \ldots \}$ belonging to the range of $-f'$. 
Each $C_k$ is a (topological) hyperplane, which admits a simple explicit parametrization 
by functions of average zero.}

\bigskip

If $1$ is the only square integer in the range of $-f'$, 
Theorem A is part of the proof of the original Ambrosetti-Prodi Theorem. 
If $f'$ crosses only the first $m$ square integers, 
the result was obtained by Ruf (\cite{R},Proposition 9), by using different techniques.
In a forthcoming paper, Burghelea, Saldanha and Tomei prove that, for arbitrary generic 
nonlinearities, the components of the critical set are
still topological hyperplanes parametrized by square integers in the range of $-f'$. The
arguments do not have the same expliciteness than those presented in this paper, and
depend strongly on special properties of infinite dimensional topology.

In a sense, this paper is a nonlinear version of oscillation theory: 
as the function $u$ varies in $X$, consider the argument of a nonzero solution of
the linearized equation $ -v'' + f'(u)v = 0,\ v(0) = 0$ at $t = \pi$. 
It turns out that this argument
has monotonicity properties similar to those of the usual argument of a Sturm-Liouville 
solution when the potential varies so that it is increased pointwise.

We are not concerned, in this paper, with the study of the image of the critical sets $C_k$. 
In the $n$-dimensional case, for $f''>0$ or $f''<0$, the image of $F(C_1)$ is studied in \cite{BZ}.
For arbitrary interactions between $f'$ and $\sigma(\Delta)$, 
the image $N:=F(C_1)$ turns out to be a codimension 1 manifold, 
which is globally parametrized by the functions of average zero. 
For functions $g$ on one side of $N$ 
there is no solution $u\in X$ for the equation $F(u)=g$, 
while on the other side there is at least one solution. 
This result was obtained by Berger-Podolak (\cite{BP}) in the original Ambrosetti-Prodi context, 
i.e., when $C=C_1$.

When $f'$ interacts with the $j$ first eigenvalues of the free Laplacian, 
the classification of the singularities in $C_k$ is still an open problem, even in the 
one dimensional case. It is well known that $C_1$ consists only of fold points.
Also, higher singularities do appear in $C_k$, $k=2,\ldots,j$,
but the Morin type of the singularities which may occur is unknown  (see \cite{R}).

\section{Statements and proofs}

A simple computation obtains the derivative of $F$ at $u$,
\[\begin{array}{cccc}DF(u):&X&\to &Y\\ &w&\mapsto&-w^{\prime\prime}+f'(u)w.
\end{array}\]
The characterization of the critical points of $F$ is then a consequence of Fredholm theory 
applied to Sturm-Liouville operators. Thus, $u\in C$ if, and only if, the kernel of $DF(u)$ is non-trivial. 
For Dirichlet boundary conditions, the spectrum is simple.
Define $v(u)(t)$ as the solution of the linearized equation
\begin{equation}\label{problin}
-[v(u)]''(t)+f'(u(t))v(u)(t)=0,\quad v(u)(0)=0,\ [v(u)]'(0)=1.
\end{equation}
The simplicity of the spectrum of $DF(u)$ guarantees that, if $v(u)(\pi)=0$, then $\ker DF(u)$ is spanned by $v(u)(t)$.
Following Pr\"ufer (\cite{CL}), let $W:X\times [0,\pi]\to \RR$ be the continuously defined argument of the planar vector 
$([v(u)]'(t),v(u)(t))$, with $W(u)(0)=0$. It follows that $u\in C$, the
critical set, if, and only if, $W(u)(\pi)=k\pi$, $k\in\NN^*:=\{1,2,3,\ldots\}$. 
We define, for all $\theta\in \RR$, the level sets
\[M_\theta:=\{u\in X;\, W(u)(\pi)=\theta\}.\]
The critical set then decomposes into
\[C_k:=M_{k\pi}=\{u\in X;\, W(u)(\pi)=k\pi,\ k\in \NN^*\}.\]
\bigskip

\noindent {\bf Theorem B}: {\it Let $f\in C^r(\RR)$, $r\geq 2$ be a function such that $f''(0)\neq 0$. Then $M_\theta$ (in particular, each $C_k$) is either empty or a $C^r$-manifold of codimension $1$ in $X$. 
If, however, $f''(0)=0$ and the two conditions below hold,
\begin{enumerate}
\item [$(a)$] the root $0$ of $f''$ is isolated,
\item [$(b)$] $f'(0)\neq -j^2$, $j\in \NN^*$,
\end{enumerate}
then the same conclusion is valid for the critical sets $C_k$}.

\begin{proof}
We calculate $DW(u)(\pi)$.  Differentiation produces
\begin{equation*}\label{eqDW(u)}
DW(u)(\pi)\cdot\varphi=\frac{[v(u)]^{\prime}(\pi).[Dv(u)\cdot\varphi](\pi)-v(u)(\pi).[Dv(u)\cdot\varphi]^{\prime}%
(\pi)}{\{v(u)(\pi
)\}^{2}+\{[v(u)]^{\prime}(\pi)\}^{2}}.
\end{equation*}
In order to obtain $[Dv(u)\cdot\varphi](\pi)$ and $[Dv(u)\cdot\varphi]^{\prime}(\pi)$, we differentiate problem $(\ref{problin})$ with respect to $u$: for all $\varphi\in X_{D}^{2}$,
\begin{equation*}
\begin{array}
[c]{l}%
-D\big([v(u)]^{\prime\prime}\big)\cdot\varphi+ [f^{\prime\prime
}(u).\varphi].v(u)+f^{\prime}(u).(Dv(u)\cdot\varphi)=0\\
(Dv(u)\cdot\varphi)(0)=0,\quad[Dv(u)\cdot\varphi]^{\prime}(0)=0,
\end{array}
\end{equation*}
Denoting $\mu=\mu(u,\varphi)=Dv(u)\cdot\varphi\in X^{2}$, we are thus led to the initial value problem 
\begin{equation*}
\begin{array}
[c]{l}%
-\mu^{\prime\prime}+f^{\prime}(u).\mu=-[f^{\prime\prime}(u).\varphi].v(u)\\
\mu(0)=0,\quad\mu^{\prime}(0)=0,
\end{array}
\end{equation*} 
which can be solved by variation of constants. 
The expressions we want to calculate are evaluations of $\mu$ and $\mu^{\prime}$ at $\pi$. We find
\[DW(u)(\pi)\cdot\varphi=\frac{-1}{\{v(u)(\pi)\}^{2}+\{[v(u)]^{\prime}(\pi
)\}^{2}}\int_{0}^{\pi}f^{\prime\prime}(u(r)).\varphi(r).[v(u)(r)]^{2}\ dr.\]
Thus $DW(u)(\pi)\equiv0$ if, and only if, $f^{\prime\prime}(u)\equiv0$, 
which may happen only if $f''(0)=0$. So, if $f''(0)\neq 0$, 
there exists $\varphi\in X_{D}^2$ such that $DW(u)(\pi)\cdot \varphi\neq 0$ and 
the result now follows from the Implicit Function Theorem. 

Suppose now $f''(0) = 0$. Since $0$ is an isolated root of $f''$ and we must have $f''(u(t)) \equiv 0$, 
it follows that $u\equiv 0$. But, in this case, denoting $c=f'(0)$, we see that $v=v(u)$ solves 
\[-v''+cv=0,\quad v(0)=0,\quad v'(0)=1\]
and $v(u)(\pi)\neq 0$ if and only if $c\neq -j^2$. Again, $C_k$ is a manifold. 
\end{proof}

Let us now consider the Ambrosetti-Prodi context $f''>0$ or $f''<0$. 
Let $p(t)\in X$ be any strictly 
positive function in $(0,\pi)$ and $V$ be the space spanned by $p(t)$. 
Decompose $X=H\oplus V$ in $L^2-$orthogonal terms.
\bigskip

\noindent {\bf Theorem C}:  {\it If $f''>0$ or $f''<0$, then, for $j \in \NN^*$,  
$C_j\neq\emptyset$ if and only if 
$-j^2$ belongs to the interior of the image of $f'$. Also, if $M_\theta\neq\emptyset$ (in 
particular, $C_k = M_{k \pi} \neq\emptyset$), then the projection 
\[\Pi:X=H\oplus V\to H.\]
is a diffeomorphism from $M_\theta$ to $H$.}
\bigskip

\begin{proof}
For each $\lambda\in \RR$ and $h\in H$ fixed, we consider the straight line 
$\ell_{h,p} = h+\lambda p\in X$. 
We will show that $\ell_{h,p}$ always intercepts each manifold $M_\theta$ once and transversally.
Uniqueness and transversality follows from 
\[DW(u)(\pi)\cdot p=\frac{-1}{\{v(u)(\pi)\}^{2}+\{[v(u)]^{\prime}(\pi
)\}^{2}}\int_{0}^{\pi}f^{\prime\prime}(u(r))p(r)[v(u)(r)]^{2}\ dr \neq 0.\]
Smoothness (and local smooth invertibility, for a fixed $\theta$) in $h$, in turn, 
follows by setting $u = h + \lambda p$ in the formula above.

Suppose that $f''<0$: in this case, $-f'$ is increasing, $-f'(-\infty)=a$ and 
$-f'(\infty)=b$. Here $a,b\in [-\infty,\infty]$. Clearly, the range $(a,b)$ of $-f'$ 
is an open set. 
From standard oscillation theory, the solutions of the three problems below 
(with initial position and velocity at 0 equal respectively to 0 and 1)
have increasing arguments,
\[ v_a'' + a v_a = 0 \quad v_\lambda'' - f'(h + \lambda p) v_\lambda = 0 \quad v_b'' + b v_b =0.\]
We will see that the asymptotics of the argument at both ends of $\ell_{h,p}$
will be given by the argument of the solutions of the leftmost and rightmost problems.

For a fixed value of the parameter $\omega$, let $v_\omega$ be the solution of the problem
\[v''+\omega v=0,\quad v(0)=0,\ v'(0)=1.\]
We denote
\[W_1=W_1(\omega,t)=W(v_\omega)(t)\quad\text{and}\quad W_2=W_2(\lambda,t)=W(h+\lambda p)(t).\]
It follows immediately that 
\[\lim_{\omega\to -\infty}W_1(\omega,\pi)=0\quad\text{and}
\quad\lim_{\omega\to\infty}W_1(\omega,\pi)=\infty.\]

We first study the behavior of $W_2(\lambda,t)$ when $\lambda\to\infty$. 
Fix $\omega\in\RR$. It is easy to prove (see \cite{CL}) that the argument function satisfies the differential equation
\[[W(u)(t)]'={\rm cos}^2 W(u)(t)-f'(u(t)){\sin}^2 W(u)(t).\]
Then

\[(W_1-W_2)'\!=\![\omega\!+\!f'(h+\lambda p)]
\sin^2\!W_2\! + \left(\frac{\cos^2\!W_1\!-\!\cos^2\!W_2}{W_1-W_2}\right)\!(W_1\!-\!W_2)
.\]
Thus $U^{\prime}+gU=H$, for
$U:=W_1-W_2$, $H=:[\omega+f^{\prime}(h+\lambda p))]\sin^2 W_2$ and  
\[g:=-\left(\frac{\cos^2W_1-\cos^2W_2%
}{W_1-W_2}\right).\]
Since $g$ is continuous and uniformly bounded in $\lambda$, we can define the integrating factor $G(t)=\exp(\int_0^t g(s)ds)$. Multiplication by this factor and integration produces
\begin{equation}\label{lemma}
\exp(G(\pi))U(\pi)=\int_0^\pi \exp(G(t))H(t)dt.
\end{equation}
If $b<\infty$, we choose $\omega=b$ and have $H(t)>0$, yielding $W_1(b)>W_2(\lambda)$ for all $\lambda$. The estimate
\[\int_{0}^{\pi}\exp (G(t))H(t)dt\leq \left( M\|\sin^2W_2\|_{L^{2}}\right) \| b-(-f^{\prime }(h+\lambda p))\|_{L^{2}},\]
where $M:=\max_{t\in [0,\pi ]}\exp (G(t))$, shows that $U(\pi)\to 0$ when $\lambda\to\infty$. 
This shows that $W_2(\lambda, \pi)$ converges increasingly to $W_1(b,\pi)$, 
for $\lambda \to \infty$, if $b < \infty$.
If $b=\infty$, we have $-f^{\prime }(h+\lambda p)-\omega >0$ in the interval $(\delta,\pi-\delta)$ if $\lambda$ is big enough. Defining $r(t):=\sin^2W_2\exp(G(t))$, we obtain
\[\exp (G(\pi))(U(\pi ))\leq 2M\delta\omega-\int_{\delta}^{\pi-\delta}[-f'(h+\lambda p)-\omega]r(t)dt<0\]
when $\lambda\to\infty$. Consequently, the integral in (\ref{lemma}) is negative, 
if $\lambda$ is sufficiently big. We conclude that $W_2(\lambda,\pi)>W_1(\omega,\pi)$. 
Since $\lim_{\omega\to\infty}W_1(\omega,\pi)=\infty$, and we again have that 
$W_2(\lambda, \pi)$ increases to $W_1(b,\pi) = \infty$. 

In a similar fashion, we obtain that $W_2(\lambda,\pi)$ converges decreasingly to 
$W_1(a,\pi)$, when 
$\lambda \to - \infty$. Thus, the line $\ell_{h,p}$ intersects the same manifolds $M_\theta$
irrespective of $h$: the values of $\theta$ for which $\ell{h,p}$  trespasses $M_\theta$ 
lie strictly between $a$ and $b$. The case $f'' >0$ is similar.
\end{proof}


\begin{thebibliography}{99}
\bibitem{AP}  A. Ambrosetti and G. Prodi: \textit{On the inversion of
some differentiable mappings with singularities between Banach spaces}, Ann.
Mat. Pura Appl. \textbf{93} (1972), 231-246.
\bibitem{BP}  M. S. Berger and E. Podolak: \textit{On the solutions of a nonlinear Dirichlet problem}, Indiana Univ. Math. J. \textbf{24} (1975), 837-846.
\bibitem{BZ} H. Bueno and A. Zumpano: \textit{Remarks about the geometry of the Ambrosetti-Prodi problem}, to appear in Nonlinear Analysis, TMA. 
\bibitem{C} V. Cafagna: {\it Global invertibility and finite solvability}, Nonlinear Funct. Anal. Lecture Notes Pure Appl. Math. Vol. 121, Ed. P. S. Milojevi\'{c}, Marcel Dekker, New York, 1990.
\bibitem{CL} E. Coddington and N. Levinson: Theory of Ordinary Differential Equations, McGraw-Hill Book Company (1955), T\ M\ H Edition, Bombay. 
\bibitem{MST}  I. Malta, N. C. Saldanha and C. Tomei: \textit{Morin singularities and global geometry in a class of ordinary differential operators}, Top. Meth. Nonlinear Anal. \textbf{10} (1997), 137-169.
\bibitem{M} H. McKean: {\it Singularities of a simple elliptic operator}, J. Differential Geometry \textbf{25} (1987), 157-165 and {\it Correction to ``Singularities of a simple elliptic operator"}, J. Differential Geometry \textbf{36} (1992), 255.
\bibitem{R} B. Ruf: {\it Singularity theory and bifurcation phenomena in differential equations}, Topological nonlinear Analysis II (Frascati, 1995), 315-395, Progr. Nonlinear Differential Equations Appl. {\bf 27}, Birkh\"{a}user Boston, Boston, 1997.
\end{thebibliography}
\end{document}